%%%%%%%%%%%%%%%%%%%%%%%%%%%%%%%%%%%%%%%%%%%%%%%%%%%%%%%%%%%%%%%%%%%%
%																   %
%  Miran Cerne, Franc Forstneric								   %
%																   %
%  EMBEDDING SOME BORDERED RIEMANN SURFACES IN THE AFFINE PLANE	   %
%																   %
%  Last edited: May 22, 2002					    			   %
%  																   %
%  Math. Res. Lett.                        		   %
%																   %
%%%%%%%%%%%%%%%%%%%%%%%%%%%%%%%%%%%%%%%%%%%%%%%%%%%%%%%%%%%%%%%%%%%%
%
%  This is the Macro file.
%
%
\scrollmode
\magnification=\magstep1
\parskip=\smallskipamount

\hoffset=1cm
\hsize=12cm

\def\demo#1:{\par\medskip\noindent\it{#1}. \rm}
\def\ni{\noindent}               % noindent
     % new indented par in horizontal mode
 % new unindented par in horizontal mode
\def\ll{\leftline}
\def\cl{\centerline}

\def\begin{\ll{}\vskip 10mm \nopagenumbers}  % beginning of the paper
\def\pn{\footline={\hss\tenrm\folio\hss}}   % pagenumbers at bottom
\def\ii#1{\itemitem{#1}}

%
%   \beginsection
%
\outer\def\beginsection#1\par{\bigskip
  \message{#1}\leftline{\bf\&#1}
  \nobreak\smallskip\vskip-\parskip\noindent}

%
%  \proclaim
%
\outer\def\proclaim#1:#2\par{\medbreak\vskip-\parskip
    \noindent{\bf#1.\enspace}{\sl#2}
  \ifdim\lastskip<\medskipamount \removelastskip\penalty55\medskip\fi}

\def\endpr{\hfill $\spadesuit$ \medskip}

%
%
%  Roman capital characters
%
%
               % roman D in math mode

               % roman T
               % roman N in math

%
%  Special letters: real, complex numbers etc.
%

\def\R{{\rm I\kern-0.2em R\kern0.2em \kern-0.2em}}
\def\N{{\rm I\kern-0.2em N\kern0.2em \kern-0.2em}}
\def\P{{\rm I\kern-0.2em P\kern0.2em \kern-0.2em}}
\def\B{{\rm I\kern-0.2em B\kern0.2em \kern-0.2em}}
\def\C{{\rm C\kern-.4em {\vrule height1.4ex width.08em depth-.04ex}\;}}
\def\CP{\C\P}
% \def\T{{\rm I\kern-0.5em T\kern0.2em \kern-0.2em}}

%
%
% Caligraphic capital characters
%
%
\def\cA{{\cal A}}
\def\cB{{\cal B}}
\def\cC{{\cal C}}

\def\cF{{\cal F}}

\def\cJ{{\cal J}}

\def\cM{{\cal M}}

\def\cO{{\cal O}}

\def\cR{{\cal R}}

\def\cT{{\cal T}}

\def\cW{{\cal W}}

%
%
%  Small Greek letters in Math mode
%
\def\a{\alpha}

\def\d{\delta}
\def\e{\epsilon}

\def\k{\kappa}

%\def\r{\rho}

%\def\o{\omega}

%\def\O{\Omega}

%
%
%  Miscellaneous
%
%
\def\bar{\overline}              % conjugate
\def\bs{\backslash}              % backslash
             % disc
        % closed disc
\def\di{\partial}                % partial derivative
\def\dibar{\bar\partial}         % di-bar derivative

%
% Abbreviations
%
\def\dim{{\rm dim}\,}                    % dimension
\def\holo{holomorphic}                   % holomorphic
                  % automorphism
                  % homomorphism
               % analytic subset
                % homeomorphism
\def\nbd{neighborhood}                   % neighborhood
                   % pseudoconvex
        % strongly psc
\def\ra{real-analytic}                   % real-analytic
               % plurisubharmonic

                    % totally real
\def\pc{polynomially convex}             % polynomially convex
          % holo convex
           % holomorphic function
\def\ss{\subset\!\subset}                % relatively compact subset
                  % support
             % C^n equivalent

\def\iff{if and only if}

\def\hvb{holomorphic vector bundle}

\def\phe{proper holomorphic embedding}

\def\dist{{\rm dist}}
\def\phm{proper\ holomorphic\ map}

                         % roman Aut in math mode
              % sypmlectic automorphisms
\def\End{{\rm End}}

\def\wt{\widetilde}
\def\wh{\widehat}

%
%  Od Mirana:  Script R s krozcem. Uporabljaj kot $\IR$.
%
\def\IR{{\cal R}\kern-2.7mm\hbox{\raise0.8mm\hbox{\tenrm\char23}}}

\begin
\nopagenumbers
\centerline{\bf EMBEDDING SOME BORDERED RIEMANN SURFACES} 
\centerline{IN THE AFFINE PLANE}
\bigskip\medskip
\cl{Miran \v Cerne and Franc Forstneri\v c}
\medskip\bigskip

\rm
\beginsection 1. Introduction.

It is a long-standing open problem whether every open 
(non-compact) Riemann surface $\cR$ admits 
a \phe\ into $\C^2$. The most general result in this direction 
was proved by Globevnik and Stens\o nes [GS]: 
{\it Every  finitely connected domain in $\C$ with no isolated 
points in the boundary embeds in $\C^2$.} Earlier constructions
of embeddings in $\C^2$ are due to Kasahara and Nishino
for the disc $U=\{z\in \C\colon |z|<1\}$ [Ste], 
Laufer [Lau] for the annuli $A=\{1<|z|<r\}$ and Alexander [Ale] 
for the punctured disc $U\bs \{0\}$. Every open Riemann surface 
embeds in $\C^3$ [GR], and every compact Riemann surface embeds 
in $\CP^3$ but most of them don't embed in $\CP^2$ [FK].

In this paper we consider the embedding problem
for {\it bordered Riemann surfaces}. The underlying
space is a compact, orientable, smooth real surfaces $\cR$ with 
boundary $b\cR=\cup_{j=1}^m C_j$ consisting of finitely many 
curves. Topologically $\cR$ is equivalent to a sphere 
with $g_{\cR}$ handles (the {\it genus} of $\cR$) 
and $m\ge 1$ discs removed. A complex structure on $\cR$ is determined 
by a real endomorphism $J$ of the tangent bundle $T\cR$ 
satisying $J^2=-Id$ (Gauss-Ahlfors-Bers; for higher dimensions 
see Newlander and Nirenberg [NN]). Two complex structures 
$J_0$ and $J_1$ are equivalent if there exists a diffeomorphism
$\phi\colon \cR\to\cR$ satisfying $d\phi\circ J_0=J_1\circ d\phi$.
The set of equivalence classes of complex structure on $\cR$,
denoted $\cM(\cR)$, is called the {\it moduli space} 
of Riemann surface structures on $\cR$. We emphasize that
we do not deal with punctured surfaces, i.e, $\cR$ has no 
isolated points in the boundary.

Our first result is that the topology of a bordered 
surface plays no role in the problem of holomorphic
embeddability in $\C^2$.

\proclaim 1.1 Theorem:  On each bordered surface $\cR$
there exists a complex structure such that the 
open Riemann surface $\IR=\cR\bs b\cR$ admits 
a \phe\ in $\C^2$.

We give an elementary proof of theorem 1.1 in section 2.
In fact we will show that %, in a natural sense to be described,
{\it there is a non-empty open set of non-equivalent complex 
structures on $\cR$ for which $\IR$ embeds into $\C^2$}
(theorem 1.5).

The precise regularity of $\cR$ up to the boundary is not
important since $\cR$ is biholomorphically
equivalent to a domain bounded by $m\ge 1$ 
real-analytic curves in a compact Riemann surface 
$\wt \cR$. The simplest way to obtain such $\wt \cR$ 
is to fill each hole of $\cR$ by attaching a disc
$D_j$ such that we identify $bD_j$ with 
a boundary curve $C_j\subset b\cR$ and extend 
the complex structure across $C_j$ using reflection. 
Another possibility is to embed $\cR$ in the (Shottky) 
{\it double} $\wh\cR$ which is obtained by gluing 
two copies of $\cR$ (with the opposite orientation) along 
$b\cR$. The genus of $\wh\cR$ equals $2g_{\cR}+m-1$. 
The details of this `doubling construction' can be found in 
[BS, pp.\ 581-2], [Spr, p.\ 217] or [SS]. Hence we shall assume 
from now on that $\cR$ is smooth up to the boundary.

We shall now give a more precise embedding theorem.
For each $k\ge 0$ we denote by $\cA^k(\cR)$ the algebra of 
all $\cC^k$-functions on $\cR$ which are \holo\ in $\IR$. 
A nonconstant function $f \in \cA(\cR)$ satisfying $|f|=1$ on 
$b\cR$ is called an {\it inner function} on $\cR$. The restriction 
of an inner function to $\IR$ is a proper \holo\ map of $\IR$ 
onto the disc $U$; conversely, any \phm\ $f\colon \IR\to U$ extends 
to an inner function on $\cR=\IR \cup b\cR$.
There is an integer $d=deg(f)\in \N$, called the 
{\it degree} (or {\it multiplicity}) of $f$, 
such that for all except finitely many points $z\in \bar U$ 
the fiber $\cR_z=f^{-1}(z)$ consists of $d$ distinct points 
in $\cR$ while the exceptional fibers consist of less than $d$ points. 

\pn

%
%  class F
%
\proclaim Definition 1: A bordered Riemann surface $\cR$ 
of genus $g_{\cR}$ and with $m$ boundary components is 
said to be of class ${\cF}$ if it admits an injective immersion 
$F=(f,g) \colon \cR\to \bar U\times \C$ which is 
\holo\ in $\IR$ and such that $f$ is an inner function 
with $deg(f) \ge  2g_\cR+m-1$.

Clearly this property is biholomorphically invariant.
The reason for introducing this class is the following 
result which is proved in section 3. 

\proclaim 1.2 Theorem: If $\cR$ is a bordered Riemann surface 
of class $\cF$ then $\IR=\cR\bs b\cR$ admits a \phe\ in $\C^2$.

\ni \it Examples: \rm 1. On each smoothly bounded domain 
$\Omega\ss\C$ with $m$ boundary components there exists 
an inner function $f$ with $deg(f)=m$ [Ahl]. 
The map $F(x)=(f(x),x) \in\C^2$ for $x\in\bar\Omega$ satisfies 
the hypothesis of theorem 1.2 and hence $\Omega$ embeds in $\C^2$.
This is the theorem of Globevnik and Stens\o nes [GS].

\ni 2. A compact Riemann surface is called {\it hyperelliptic} 
if it is the normalization of a complex curve 
in $\CP^2$ given by $w^2=\Pi_{j=1}^k (z-z_j)$ for some choice 
of points $z_1,\ldots,z_k\in \C$ [FK]. We shall call a bordered 
Riemann surface $\cR$ hyperelliptic if its double $\wh\cR$
is hyperelliptic. Such $\cR$ has either one or two boundary 
components and it admits a pair of inner functions $(f,g)$ which embed 
$\IR$ in the polydisc $U^2$ such that $b\cR$ is mapped to 
the torus $(bU)^2$; moreover, one of the two functions has 
degree $2g_{\cR}+m$ and the other one has degree $2$
(see [Ru1] and sect.\ 2 in [Gou]). 
Thus $\cR$ is of class $\cF$ and we get

\proclaim 1.3 Corollary: If $\cR$ is a hyperelliptic bordered 
Riemann surface then $\IR$ admits a \phe\ in $\C^2$. 
In particular, each torus with one hole embeds 
properly holomorphically into $\C^2$.

Indeed it is shown in [Gou] that the double $\wh\cR$ of a 
hyperelliptic bordered Riemann surface $\cR$ can be represented
as a curve in $\CP^2$ given by the equation
$$ 
	y^2=\prod_{j=1}^{\hat g+1} (x-\a_j)(1-\bar{\a_j}x) \eqno(1.1)
$$
for some choice of distinct points $\a_j\in U$ ($1\le j\le \hat g+1$), 
where $\hat g=2g_{\cR}+m-1$ is the genus of $\wh\cR$, such 
that $\cR=\{(x,y)\in \wh\cR\colon |x|\le 1\}$. The 
pair of inner functions on $\cR$ given by
$f=y/\prod_{j=1}^{\hat g+1} (1-\bar{\a_j}x)$, $g=x$,
provides an embedding $F=(f,g)\colon \cR \to (\bar U)^2$.
Clearly $g$ has multiplicity $2$ on $\cR$. From the 
relation $f^2=\prod_{j=1}^{\hat g+1}(g-\a_j)/(1-\bar{\a_j}g)$ 
which follows from (1.1) we see that $f$ has multiplicity 
$\hat g+1=2g_{\cR}+m$ and hence theorem 1.2 applies. 
Sikorav [Sik] gave a slightly different proof of 
corollary 1.3 for tori with one hole (unpublished);
these are all hyperelliptic.

\demo Sketch of proof of theorem 1.2: 
Let $P\ss \C^2$ be a polydisc. According to [Glo] 
(see also [\v CG] and [Stn]) there exist Fatou-Bieberbach domains 
$\Omega\subset\C^2$ such that $b\Omega\cap P$ is an arbitrarily 
small perturbation of the cylinder $(bU\times \C)\cap P$
(proposition 2.1 below). Such domains $\Omega$ can be constructed 
using sequences of compositions of shears in coordinate directions 
on $\C^2$. Let $\phi\colon\Omega\to\C^2$ be a biholomorphic map 
onto $\C^2$. If $F=(f,g)\colon\cR\to \bar U\times\C$ is as in 
theorem 1.2, we choose the polydisc $P$ large enough to contain 
$F(\cR)$ and solve a Riemann-Hilbert boundary value problem to 
find a new \holo\ embedding 
$\wt F=(\wt f, g)\colon \cR\to\bar \Omega \cap P$ 
such that $\wt F(b\cR) \subset b\Omega$. The map 
$G=\phi\circ \wt F\colon \IR\to\C^2$ is then a
\phe\ of $\IR$ to $\C^2$. The details are carried out
in section 3. This approach was used in [\v CG] for 
finitely connected planar domains to provide a different 
proof of the embedding theorem of Globevnik and Stens\o nes [GS].
The proof of theorem 1.1 (section 2) is similar but 
more elementary.
\endpr

The proof of theorem 1.1 shows that each bordered 
surface $\cR$ carries a complex structure $J$ such that 
$(\cR,J)$ is of class ${\cF}$. We will show that
the set of such complex structures on $\cR$ is open.
To be precise, fix a number $\a\in (0,1)$
and denote by $\End_{\R}^\a (T\cR)$ the set of all 
endomorphisms of $T\cR$ which are H\"older continuous of 
class $\cC^\a$, endowed with the $\cC^{\a}$ topology. Let 
$$ \cJ_\cR^\a =\{ J\in  \End_{\R}^\a (T\cR)\colon J^2=-Id\}. $$
% denote the set of all complex structures of class $\cC^\a$ on $\cR$.

\proclaim 1.4 Theorem: Let $\cR$ be a smooth bordered 
surface and $0< \a<1$. The set $\cF^\a_\cR$, consisting of all 
$J\in \cJ_\cR^\a$ such that the Riemann surface $(\cR,J)$ 
is of class $\cF$, is open in $\cJ_{\cR}^\a$.

Theorem 1.4 is proved in section 4. The main point is that
inner functions on $(\cR,J)$ with multiplicity
at least $2g_{\cR}+m-1$ are stable under small perturbations
of the complex structure $J$ (proposition 4.1).

For each $\a\in (0,1)$ we can realize the moduli space 
$\cM(\cR)$ as the quotient $\cM(\cR)=\cJ_{\cR}^\a / \!\!\sim$ where 
$J_0,J_1 \in  \cJ_{\cR}^\a$ satisfy $J_0\sim J_1$ \iff\ there exists
a diffeomorphism $\phi\colon \cR\to\cR$ of class $\cC^{1,\a}$ 
with $d\phi\circ J_0=J_1\circ d\phi$.
Let $\pi\colon 	\cJ_{\cR}^\a \to \cM(\cR)$ denote the quotient
projection. We endow $\cM(\cR)$ with the quotient topology.
The set $\cF(\cR)=\pi(\cF^\a_\cR) \subset \cM(\cR)$ consists
of moduli of Riemann surface structures of class $\cF$;
this makes sense since the property $\cF$ is
biholomorphically invariant and therefore 
$\cF^\a_\cR =\pi^{-1}(\cF(\cR))$. We can now summarize
the above results as follows.

\proclaim 1.5 Main Theorem: Let $\cM(\cR)$ denote
the moduli space of Riemann surface structures on a 
smooth bordered surface $\cR$. 
The set $\cF(\cR) \subset \cM(\cR)$, 
consisting of structures of class $\cF$, is nonempty and 
open in $\cM(\cR)$. For each $[J]\in \cF(\cR)$ the open 
Riemann surface $(\IR,J)$ admits a \phe\ in $\C^2$.
If $\cR$ is a finitely connected domain in $\C$ 
then $\cF(\cR)=\cM(\cR)$.

The last statement above is the theorem of [GS]. The main question 
is whether $\cF(\cR) =\cM(\cR)$ for every bordered surface $\cR$. 
The discussion in section 2 seems to indicate a negative answer; 
see proposition 2.2 and the remark following its proof.

\demo Comments regarding class $\cF$:
It is proved in [Ahl, pp.\ 124-126] that on every bordered 
Riemann surface $\cR$ of genus $g_{\cal R}$ with 
$m$ boundary components there is an inner 
function $f$ with multiplicity $2g_{\cal R} + m$ 
(although the so-called Ahlfors functions may have 
smaller multiplicity). A generic choice of 
$g\in A^1(\cR)$ gives an immersion 
$F=(f,g)\colon \cR\to \bar U\times \C$ with at most 
finitely many double points (normal crossings). 
The main difficulty is to find $g$ such that $F=(f,g)$ 
is injective on $\cR$. We do not know whether 
such $g$ always exists as Oka's principle does not apply 
in this situation (proposition 2.2).
\endpr

It would be of interest to relax the condition in 
theorem 1.2 that one of the components be an inner 
function. In this direction we pose the following

\medskip \ni\bf Problem: \sl 
Let $\cR$ be a bordered Riemann surface
and let $F=(f,g)\colon \cR\to\C^2$ be a
holomorphic embedding whose image 
$F(\cR)$ is \pc\ in $\C^2$. Does $\IR$ 
embed (properly holomorphically) in $\C^2$ ? 
\rm \medskip

A possible approach would be to use sequences of automorphisms 
of $\C^2$ to carry the boundary points of $F(b\cR)$ towards
infinity. 

The problem of holomorphic embeddability of a bordered
Riemann surface $\cR$ in $\C^2$ is related to the question 
whether certain algebras of holomorphic functions on 
$\cR$ are doubly generated. If $F=(f,g)\colon \cR\to\C^2$
is an $\cA^k$-embedding whose image $F(\cR)$ is \pc\ then 
$f$ and $g$ generate a dense subalgebra of $\cA^k(\cR)$; 
in such case we say that the algebra is {\it doubly generated}.
Conversely, if $f$ and $g$ generate a dense subalgebra
in $\cA^k(\cR)$ then $F=(f,g)\colon \cR\to\C^2$ 
is an injective immersion (not necessarily proper). 
The question whether $\cA(\cR)$ is always doubly generated is 
to our knowledge still open. In 1978 Tsanov [Tsa] proved that 
for any bordered surface $\cR$ there is a non-empty open set in 
the Teichm\"uller space $\cT(\cR)$ (or in the reduced
Teichm\"uller space $\cT^\#(\cR)$) consisting of Riemann surfaces 
for which the algebra $\cA^0(\cR)$ is doubly generated. In view of 
the above remark this also follows from theorem 1.5.

Higher dimensional analogues of open Riemann surfaces
are {\it Stein manifolds}. A complex manifold is Stein 
if it admits a \phe\ in some complex Euclidean space $\C^N$
[GR, H\"or]. It is known that every Stein manifold of complex 
dimension $n\ge 2$ admits a \phe\ into $\C^N$ with 
$N=[3n/2]+1$ and this $N$ is in general the smallest possible 
[EG, Sch]. For $n=1$ this would give $N=2$; however, the proof, 
which is based on the `removal of singularities' method, 
does not apply for $n=1$ since its main ingredient 
breaks down. This crucial ingredient is the {\it homotopy principle}, 
also called the {\it Oka-Grauert principle}, for sections of 
\hvb s over Stein manifolds which avoid certain complex 
analytic subvariety of the total space [Gra, Gro2, FP1, FP2, FP3]. 
When $n=1$, the subset to be avoided is a complex hypersurface 
and the Oka-Grauert principle fails in general (proposition 2.2).

\demo Acknowledgements: 
We wish to thank Josip Globevnik and Berit Stens\o nes for 
stimulating discussions on this topic. This research was 
supported in part by the Ministry of 
Science of the Republic of Slovenia.

\beginsection 2. Proof of theorem 1.1.

In this section we outline an approach to construct
embeddings of bordered Riemann surfaces $\cR$ in the tube 
$\bar U\times \C$ by the `elimination of singularities' 
method which has been used successfuly in higher dimensions
[EG, Sch]. Although we cannot prove the existence
of such an embedding for every complex structure on a
given surface $\cR$, we obtain an elementary proof of 
theorem 1.1 based on the following result of Globevnik [Glo]. 
We quote the version proved in [\v CG]. As before, $U$ denotes 
the unit disc in $\C$ and $rU=\{|z|<r\}$. We denote the coordinates 
on $\C^2$ by $(z,w)$. Fix $r>0$ and let 
$P=(2U)\times (rU) \subset\C^2$.

%
%  Globevnik's theorem
%
\proclaim 2.1 Proposition: {\rm ([\v CG], Lemma 2.1)} 
There exist arbitrarily small $\cC^3$-pertur\-bations 
$S$ of $bU\times \C$ such that, if $\Omega$ is the 
domain in $P$ bounded by $S\cap P$ and by 
$|w|=r$, there is an injective \holo\ map 
$\phi\colon\Omega\to\C^2$ such that 
$|\phi(z_n,w_n)|\to +\infty$
whenever $(z_n,w_n)\in \Omega$ %($n\in \N$) 
and $\lim_{n\to\infty} \dist\{(z_n,w_n),S\} =0$.

\demo Remark: In fact there exist Fatou-Bieberbach domains
$\Omega\subset\C^2$ with smooth boundary  such that 
$b\Omega\cap P$ is a small perturbation of 
$(bU\times \C)\cap P$ (Stens\o nes [Stn]; see also
Globevnik [Glo] for the $\cC^1$ version). The weaker result 
quoted above will suffice for our purposes.

%
%  Proof of 1.1.
%
\demo Proof of theorem 1.1:	
We may assume that $\cR$ is a compact domain with smooth 
\ra\ boundary in a Riemann surface $\wt \cR$ (sect.\ 1). 
We denote by $\cO(\cR)$ the algebra of functions 
$f$ holomorphic in a \nbd\ $V_f$ of $\cR$ in $\wt\cR$.
By Ahlfors [Ahl] there is an inner function $f$ on $\cR$. 
Let $Z=\{z_1,z_2,\ldots,z_p\} \subset \bar U$ be the set of critical 
values of $f$. By the Hopf lemma we have $df_x\ne 0$ for 
$x\in b\cR$ and hence $Z$ is contained in the open disc $U$. 
Choose a function $g_1 \in \cO(\cR)$ such that

\item{(a)} $g_1$ separates points on the
(finitely many) fibers $\cR_z$ for $z\in Z$, and
\item{(b)} $dg_1\ne 0$ at each point $x\in f^{-1}(Z)\subset \cR$. 

Clearly $g_1$ will separate points on all except perhaps
finitely many fibers $\cR_z$ ($z \in \bar U$). Condition 
(b) insures that $F_1=(f,g_1) \colon \cR\to \bar U\times \C$ 
is an immersion. A generic choice of $g_1$ also insures that 
$F_1$ has only finitely many double points in $\IR$ and no 
double point on $b\cR$. Now choose $g_2 \in \cO(\cR)$ 
which vanishes to second order at each point of $f^{-1}(Z)$ 
and such that the pair $(g_1,g_2)$ separates 
points on all fibers $\cR_z$ for $z\in \bar U$; clearly
such $g_2$ exists since we must satisfy the separation
condition only at finitely many points.

We wish to find $g\in \cO(\cR)$ such that 
$F=(f,g)\colon \cR\to\bar U\times \C$ is an
embedding. As in [EG] and [Sch] we seek $g$ in the form
$$ 
   g(x) = g_1(x) + \alpha(f(x))\, g_2(x), \quad x\in \cR, \eqno(2.1) 
$$
where $\alpha\colon \bar U\to\C$ is a holomorphic function
to be selected. Since $g_2$ vanishes to second order
at each point $x\in f^{-1}(Z)$, we have $g(x)=g_1(x)$
and $dg_x=(dg_1)_x$ at such points. Thus for any choice of 
$\a$ the associated map $F=(f,g) \colon\cR\to \bar U\times \C$ 
is a \holo\ immersion which is injective in a \nbd\ of
$f^{-1}(Z)$.

Our goal is to choose $\a$ such that $F$ is injective 
globally on $\cR$. To formulate the relevant condition 
on $\a$ we fix a point $z\in \bar U$ and write
$f^{-1}(z)=\{x_1,\ldots,x_d\}$ (distinct points!), 
where $d=deg(f)$ for all except finitely many $z\in \bar U$. 
Denote by $\Sigma_z \subset \C$ the (finite) set of
solutions of the equations
$$   g_1(x_i)+ag_2(x_i)= g_1(x_j)+ag_2(x_j),\qquad
     1\le i<j\le d.                               
$$
Equivalently, a number $a\in \C$ belongs to $\Sigma_z$
if it solves the equation
$$ 
   a \bigl( g_2(x_j)-g_2(x_i) \bigr) = g_1(x_i)-g_1(x_j) 
$$
for some $i\ne j$. By the choice of $g_1$ and $g_2$
at least one of the differences above 
is nonzero for each pair of indices $i,j$ and hence each 
equation has either one solution or no solutions. The set
$\Sigma=\cup_{z\in\bar U} \{z\} \times \Sigma_z \subset \bar U\times \C$                        
is a closed one-dimensional complex analytic subset 
of $\bar U\times \C$. The function $g$ determined by $\a$ 
according to (2.1) separates the points on all fibers of 
$f$ \iff\ the graph of $\a$ avoids $\Sigma$, that is, if
$\a(z) \notin \Sigma_z$ for all $z\in \bar U$.

Choose a simple smooth arc $C\subset U$ containing the 
set $Z$ of critical values of $f$. 
By dimension reasons there is a smooth function 
$\a_0\colon C\to \C$ whose graph over $C$ 
avoids $\Sigma$. We can approximate $\a_0$ uniformly on 
$C$ by \holo\ polynomials $\a$. If the approximation is 
sufficiently close, the graph of $\a$ will avoid $\Sigma$ 
over an open set $V\ss U$ containing $C$. If $g$ 
is the corresponding function (2.1) then $F=(f,g)$ is a 
proper \holo\ embedding of $\cR_V=f^{-1}(V)\subset \cR$ 
to $V \times \C$.

Choose a simply connected closed domain $D_0 \ss V$ with 
\ra\ boundary and containing the arc $C$ in its interior. 
There are a domain $V_1\subset V$ containing $D_0$ and an 
injective holomorphic map $\sigma\colon V_1\to \C$ which 
maps $D_0$ conformally onto $\bar U$. The map
$F'=(\sigma\circ f, g) \colon f^{-1}(V_1) \to \C^2$ 
is a \holo\ embedding which maps the closed domain
$\cR_0=f^{-1}(D_0) \subset \IR$ to $\bar U\times \C$
and it maps $b\cR_0$ to $bU\times \C$.

Choose a number $r> \sup\{ |g(x)|\colon x\in \cR\}$. 
Let $\Omega$ be as in proposition 2.1 such that
$b\Omega\cap (2U\times rU)$ is a small 
$\cC^3$-perturbation of the cylinder $bU \times rU$. 
If the approximation is sufficiently good then $b\Omega$ 
intersects the image of $F'$ transversely and the set 
$$
   \cR'=\{x\in  f^{-1}(V_1) \colon F'(x)\in \Omega\} \ss \IR
$$ 
is a domain in $\IR$ with $\cC^3$-boundary which is a small 
$\cC^3$-perturbation of $b\cR_0$. If $\phi\colon\Omega \to\C^2$
is chosen as in proposition 2.1 then the map
$G= \phi \circ F' \colon \cR' \to \C^2$ is a \phe\ of
$\cR'$ to $\C^2$.

To conclude the proof of theorem 1.1 it remains to show 
that $\bar {\cR'}$ is diffeomorphic to $\cR$. 
This can be seen as follows. Since $D_0$ is a closed 
simply connected domain with smooth boundary in $U$, 
there is a smooth function $\rho\colon\bar U\to \R$ 
such that $D_0=\{\rho\le 0\}$, $\rho=1$ on $bU$,
and $\rho$ has no critical points in $\bar U\bs C$ 
(hence $0<\rho(x)<1$ for $x\in U\bs D_0$). 
For $0\le t\le 1$ set 
$D_t = \{x \in \bar U\colon \rho(x) \le t\}$;
thus $D_0$ is the given set and $D_1=\bar U$.
Since $f$ has no critical values in $\bar U\bs C$, 
the function $\rho\circ f \colon \cR\to\R$ has no 
critical points in $\cR\bs f^{-1}(C)$.
By Morse theory the set 
$\cR_t=\{x\in  \cR\colon f(x)\in  D_t \} = 
\{x\in \cR \colon \rho( f(x))\le t\}$ 
is diffeomorphic to $\cR=\cR_1$ for each $t\in [0,1]$.
The same is true for $\bar \cR'$ which is a small 
$\cC^3$-perturbation of $\cR_0$. 
This completes the proof of theorem 1.1.
\endpr

\demo Remarks: 1. In the above proof we could use
the Fatou-Bieberbach domains constructed in [Glo]
whose boundaries inside a polydisc are small 
$\cC^1$-perturbations of the cylinder $bU\times U$; 
we shall need more smoothness in the proof of theorem 1.2.

\ni 2. To embed $\cR$ holomorphically into $\bar U\times \C$ 
using this scheme we would have to find a function 
$\a\in \cO(\bar U)$ whose graph avoids the complex
curve $\Sigma \subset\bar U\times \C$ constructed above. 
The fiber $\Sigma_z$ over most points $z\in U$ consists
of $d\choose 2$ points, where $d=deg(f)$. In the special
case when $d=2$ we have ${d\choose 2}=1$ and hence
the Oka-principle [Gra, Gro2, FP1, FP2, FP3] applies to 
sections of $(\bar U\times \C)\bs \Sigma$, so we obtain
a desired \holo\ function $\a$ whose graph over $\bar U$ 
avoids $\Sigma$. Note that an inner function 
$f$ of degree $d=2$ exists \iff\ $\cR$ is hyperelliptic
(since we obtain by reflection a degree two meromorphic 
function on the double $\wh \cR$ which implies that 
$\wh \cR$ is hyperelliptic). On the other hand, when $d\ge3$
the generic fiber of $\Sigma$ contains at least three points 
and hence its complement is Kobayashi hyperbolic, so 
the Oka principle does not apply.

The following result shows that there exist 
complex curves $\Sigma \subset U\times \C$ 
(not necessarily arising from our construction) 
which cannot be avoided by holomorphic graphs. 
On the other hand it is always possible to avoid 
such a curve by graph of a smooth function; hence 
{\it the Oka-Grauert principle fails for sections of
$(U\times \C) \bs \Sigma$}.

\proclaim 2.2 Proposition:
There exists a closed one-dimensional complex
subvariety $\Sigma \subset U\times \C$ which does not
contain any line $\{z\}\times \C$ and which has
a nontrivial intersection with the graph of any 
\holo\ function on $U$.

\demo Proof: Denote the coordinates on $\C^2$ by $(z,w)$.
Let $\Sigma_k \subset U\times \C$ be the union of the
following complex curves, intersected with $U\times \C$:
$$
   zw=1, \quad w=1, \quad w=jz\ \ (0\le j\le k).
$$
Assume that for each $k\in \N$ there is a \holo\ function
$\a_k$ on $U$ whose graph avoids $\Sigma_k$.
Then $\a_k$ omits the values $0$ and $1$ and hence
the sequence $\{\a_k\}_{k\in\N}$ is a normal family on $U$.
Passing to a subsequence we may assume that
$\a_k$ converges, uniformly on compacts in $U$,
to a \holo\ function $\a\colon U\to \C$
or to $\a=\infty$. 

Consider the first case. Choose numbers $0<r<1$ 
and $k_0\in \N$ such that
$$ k_0 r> \max_{|z|=r}\a(z) . \eqno(2.2) $$
For each $k\ge k_0$ the winding number of the 
function $h_k(z)= kz-\a(z)$ on the circle $|z|=r$ equals
to that of $kz$ which is one. This means that $h_k$ has a zero 
on the disc $rU=\{|z|<r\}$, i.e., the graph of $\a$
intersects the line $w=kz$ and hence $\Sigma_k$.
The same argument holds for any function satisfying (2.2).
Since for large $k\in \N$ the function $\a_{k}$ is close
to $\a$ on $r\bar U$, its graph also intersects $\Sigma_k$,
a contradiction.

In the second case when $\a_k\to \infty$
we can apply a similar argument in $U\times \bar{\C}$,
where $\bar{\C} = \C\cup \{\infty\}$ is the Riemann sphere,
to show that for all sufficiently large $k$ the graph of $\a_k$
intersects the hyperbola $zw=1$ and hence $\Sigma_k$,
a contradiction. This proves proposition 2.2.
\endpr

\ni\it Remark. \rm  
The above approach to construct a function $g$ separating 
points on the fibers of a given inner function $f$ is not quite 
as ad-hoc as it may seem. Denote by $\cO$ the sheaf of germs 
of holomorphic functions on $\IR$. According to Grauert [Gra]
the push-forward $f_*\cO$ is a coherent analytic sheaf 
of $\cO_U$-modules over the disc $U=f(\IR)$. 
For each open set $V\subset U$ we may view \holo\ 
functions on $f^{-1}(V)\subset\IR$ as sections of $f_*\cO$
over $V$. By Cartan's Theorem A [GR] the sheaf $f_*\cO$ is finitely
generated over each compact $K\ss U$, meaning that there exist 
functions $g_1,\ldots,g_n\in \cO(\IR)$ such that any  
$g\in \cO(\IR)$ may be written in the form
$g(x)=\sum_{j=1}^n \a_j(f(x))\cdotp g_j(x)$ ($x\in f^{-1}(K)$)
for some \holo\ functions $\a_1,\ldots,\a_n$ defined in a \nbd\ of $K$.
Now $g$ separates points on the fibers $f^{-1}(z)$ for $z\in K$ 
\iff\ the graph of $\a=(\a_1,\ldots,\a_n)\colon K\to\C^n$ 
avoids a complex hypersurface $\Sigma\subset U\times \C^n$
constructed as above. Proposition 2.2 indicates that this may 
not be possible in general (although we do not have a 
specific counterexample).

\beginsection 3. Holomorphic perturbations of 
bordered Riemann surfaces.

In this section we prove theorem 1.2. 
Let $P=(2U)\times U \subset \C^2$. 
For any sufficiently small perturbation $S$ of the 
cylinder $bU\times \C$ we denote by $\Omega_S$ the 
connected domain in $P$ bounded by $S$ and containing 
the origin. 

%
% Explain $\cA^k(b\cR)$, $\cA^{k,\a}(b\cR)$.
%

\proclaim 3.1 Proposition: 
Let $\cR$ be a bordered Riemann surface of genus 
$g_{\cal R}$ bounded by $m_{\cal R}$ smooth curves
and let $F_0=(f_0,g_0)\colon \cR \to \bar U\times U$ be 
a map of class $\cA^2(\cR)$ such that $f_0$ is an inner 
function on $\cR$ with $deg(f_0)\ge 2g_{\cal R} + m_{\cal R}-1$.
For any sufficiently small $\cC^3$-perturbation $S$ 
of $S_0=bU\times \C$ there is a function $f\in \cA^1(\cR)$ which 
is $\cC^1$-close to $f_0$ such that the map 
$F=(f,g_0) \colon \cR\to \C^2$ satisfies $F(\IR)\subset \Omega_S$ 
and  $F(b\cR)\subset S$.

We emphasize that the maps $F$ and $F_0$ only 
differ in the first component. If $F_0$ is an embedding,
it follows that $F$ is also an embedding provided 
that $S$ is sufficiently $\cC^3$-close to $bU\times U$.

Assuming proposition 3.1 for a moment we now prove theorem 1.2.

\demo Proof of theorem 1.2: 
Let $F_0=(f_0,g_0)\colon\cR\to \bar U\times \C$ satisfy the
hypothesis of theorem 1.2. We may assume that $\cR$ is
a domain with smooth \ra\ boundary in a larger Riemann surface 
$\wt\cR$ (section 1). Since $f_0$ maps $b\cR$ to the circle
$bU$, it extends by reflection to a \holo\ function in a \nbd\ 
of $\cR$. Furthermore we can approximate $g_0$ 
in the $\cC^1(\cR)$-sense by a function (still denoted $g_0$)
which is \holo\ in a \nbd\ of $\cR$. If the approximation is 
sufficiently close, the new map $F_0$ is an embedding 
in a \nbd\ of $\cR$. We may assume that $||g_0||_{\cR} <1$. 

Choose a hypersurface $S$ close to $S_0=bU \times U$ and the 
associated  domain $\Omega=\Omega_S\subset P$ as in proposition 2.1, 
with the corresponding injective holomorphic map 
$\phi\colon\Omega\to\C^2$ which maps sequences in $\Omega_S$ 
converging to $S$ to sequences going to infinity. Let 
$F=(f,g_0)$ be a map furnished by proposition 3.1 which is
$\cC^1$-close to $F_0$. If the approximation is sufficiently 
close then $F$ is a holomorphic embedding of $\IR$ in $\Omega$
which maps $b\cR$ to $b\Omega$. The map 
$G=\phi\circ F \colon \IR\to\C^2$ is then 
a proper holomorphic embedding of $\IR$ to $\C^2$.
This proves theorem 1.2 granted that proposition 
3.1 is correct.
\endpr

In the proof of proposition 3.1 we shall use some results about 
the linear Riemann-Hilbert problem on bordered Riemann surfaces. 
Fix a number $0<\alpha <1$. 
Denote by ${\cal A}^{1,\a}(b\cR)$  the Banach algebra 
of $\cC^{1,\a}$-functions on $b\cR$ which extend holomorphically 
to $\IR$. (The $\cC^{1,\a}$-norm on $b\cR$ can be defined 
by choosing a smooth parametrization of each curve 
$C\subset b\cR$ by the circle $S^1$ and pulling back 
functions on $C$ to functions on $S^1$.)

Given functions $a\colon b{\cal R}\rightarrow \C\bs \{0\}$ 
and $c\colon b\cR\to\R$ of class $\cC^{1,\alpha}$,
the corresponding {\it Riemann-Hilbert problem} is to find 
$k\in\cA^{1,\alpha}(b\cR)$ such that 
$$
   {\rm Re} \bigl(\, \overline {a(x)}\cdotp k(x)\bigr) =c(x),
    \qquad x \in b{\cal R}.  							\eqno(3.1)
$$
The existence of solutions depends on the {\it index} 
$\k(a)$ which is defined as the sum of the winding numbers 
of $a$ over all $m_{\cR}$ boundary components of ${\cal R}$
(the corresponding {\it Maslov index} is $2\k(a)$).
Here we equip ${\cal R}$ with the usual orientation induced by the 
complex structure and we orient the boundary $b\cR$ coherently.
Note that when $a$ is an inner function on $\cR$ we have
$\k(a)=deg(a)$. The following is a part of the theorem 
from [Kop, p.\ 30]; it corresponds to the case when $\nu=0$ 
(since we are dealing with functions) and the trivial divisor 
$\d$ with degree $n_\d=0$. Notice also that, in [Kop], the 
surface has $m+1$ holes.  

%
%  Linear R-H problem
%
\proclaim Theorem: {\rm (Koppelman [Kop])} 
Let ${\cal R}$ be a bordered Riemann surface of genus
$g_{\cal R}$ with $m_{\cal R}$ boundary components and 
let $a$ be a complex-valued H\"older continuous function 
on $b\cR$ without zeros. If $\k(a) \ge 2g_{\cal R} + m_{\cal R}-1$
then the Riemann-Hilbert boundary value problem (3.1)
is solvable for all H\"older continuous functions $c$ 
on $b{\cal R}$ and the corresponding homogeneous problem ($c=0$) 
has $2\k(a) - (2g_{\cal R} + m_{\cal R}-2)$ linearly 
independent solutions.

\demo Remarks: 1. Even though the theorem in [Kop] is stated for 
the $\cC^{\alpha}$-case (the functions $a$ and $c$ 
are assumed to be of class $\cC^{\alpha}$ and the solutions belong 
to ${\cal A}^{\alpha}(b\cR)$), the proof carries 
over to the case stated here. 

\ni 2. This is essentially a result concerning the solutions 
of the operator $L=\dibar$ acting on sections 
of the trivial line bundle $E=\cR\times \C\to \cR$, with  
the Riemann-Hilbert boundary conditions described above. 
The operator $L$ is elliptic and hence Fredholm, 
with the (real) {\it index}  
$$ {\rm Ind}(L)= 2\k(a) - (2g_{\cal R} + m_{\cal R}-2). $$
Koppelman's theorem asserts that, when  
$\k(a) \ge 2g_{\cal R} + m_{\cal R}-1$, $L$
is surjective and $\dim_{\R}(\ker L)={\rm Ind}(L)$. 
For an extension to more general $\dibar$-type operators we refer 
to [Gro1] and [HLS]. 

\ni 3. There is a connection (by doubling of $\cR$)  between
Koppelman's theorem and the Riemann-Roch theorem [HLS]. 
The result may be viewed as a special case of the Atiyah-Singer 
index theorem [AS] which expresses the index of any elliptic 
linear differential operator, acting on sections of a complex 
vector bundle $E\to X$ over a compact manifold $X$, in terms of 
the Chern class of the bundle and the cohomology class 
in $H^*(X,\C)$ determined by the principal symbol of $L$. 
% In particular, the index 
% (and the dimension of the kernel) 
% is independent of the choice of a space of sections on 
% which $L$ operates.
\endpr

\demo Proof of proposition 3.1:
Let $F_0=(f_0,g_0)$ be as in the proposition. 
Assume that $||g_0||_{\cR}<1$ so that 
$F_0(\cR) \subset P=(2U)\times U$. Denote the coordinates 
on $\C^2$ by $(z,w)$. Set $\rho_0(z,w)=|z|^2-1$.
Then $\{\rho_0=0\}\cap P= bU\times U$ and 
$\rho_0(f_0,g_0)=0$. For any function $\rho\in \cC^3(P)$ 
sufficiently close to $\rho_0$ the set 
$S_\rho=\{\rho=0\}\cap P$ is a $\cC^3$-hypersurface 
close to $bU\times U$ and vice versa, any small 
$\cC^3$-perturbation of $bU\times \C$ within $P$ 
equals $S_\rho$ for some $\rho \in \cC^3(P)$ close 
to $\rho_0$. 

To solve the problem it suffices to find for each 
$\rho\in\cC^3(P)$ close to $\rho_0$ a function 
$f=f_\rho\in\cA^{1,\a}(b\cR)$ close to $f_0$ such that 
$\rho(f(x),g_0(x))=0$ for all $x\in b\cR$. 
Such $f$ extends from $b\cR$ to a function 
$f\in \cA^{1,\a}(\cR)$. The corresponding map 
$F=(f,g_0) \colon \cR\to P$ takes $b\cR$ to 
$S_\rho=\{\rho=0\}$ and it maps the interior 
$\IR$ to the domain $\Omega_\rho\subset P$ bounded by 
$S_\rho$ and by $|w|=1$ (for the last statement we 
need $f$ to be $\cC^1$-close to $f_0$; see [\v CG] 
for the details of this argument).

To find such $f=f_\rho$ we shall apply the implicit mapping theorem 
in Banach spaces. Let $D=\{f\in \cA^{1,\a}(b\cR) \colon ||f||<2\}$, 
where $||.||$ is the norm on $\cA^{1,\a}(b\cR)$. We define a 
Banach space operator
$$ \Phi\colon \cC^3(P)\times D\to \cC^{1,\a}_{\R} (b\cR),\quad
   \Phi(\rho,f)(x) =\rho(f(x),g_0(x)) \quad (x\in b\cR).
$$
We claim that $\Phi$ is of class $\cC^1$.
Clearly $\Phi$ is linear with respect to $\rho$ and hence 
$D_\rho\Phi(\rho,f)(\tau)=\tau(f,g_0)$.
Moreover, lemma 5.1 in [HT] implies that for each 
fixed $\rho\in\cC^{3}(P)$ the mapping 
$f\in \cA^{1,\a}(b\cR) \to  \Phi(\rho, f) \in \cC^{1,\a}_{\R}(b\cR)$ 
is of class $\cC^1$ and its partial derivative on $f$ equals
$$  
    D_f\Phi(\rho,f)(k)(x)= {\rm Re} 
    \bigl( 2\di_z\rho(f(x),g_0(x))\cdotp k(x)\bigr), \quad x\in b\cR.
$$
It is easily seen that $\Phi$ and its first order partial 
derivatives are continuous with respect to both variables $(\rho,f)$ 
and hence $\Phi$ is of class $\cC^1$ [Lan]. Writing 
$a_{\rho,f}(x)= 2 \di_{\bar z} \rho(f(x),g_0(x))$ we have 
$$ 
	D_f\Phi(\rho,f)(k)(x)= 
	{\rm Re} \bigl( \, \bar{a_{\rho,f}(x)} \cdotp k(x) \bigr),
    \quad k \in \cA^{1,\a}(b\cR),\ x\in b\cR. 
$$
For $\rho_0=|z|^2-1$ we have $\di_{\bar z}\rho_0=z$ 
and hence $a_{\rho_0,f_0}=2f_0$ which is nowhere vanishing on 
$b\cR$. Its index equals   
$$ 
	\kappa(a_{\rho_0,f_0})= deg(f_0) \ge 2g_{\cR}+m_{\cR}-1.
$$
Hence the theorem quoted above applies and shows that the linear operator 
$A= D_f\Phi(\rho_0,f_0) \colon \cA^{1,\a}(b\cR)\to \cC^{1,\a}_{\R}(b\cR)$
is surjective, with kernel of dimension 
$$ \dim_\R(\ker A) = 2 deg(f_0) - (2g_{\cal R}+m_{\cal R}-2). 
$$
Since each finite dimensional subspace in a Banach space is 
complemented [Ru2, Lemma 4.21], there is a closed subspace $\cB$ in 
$\cA^{1,\a}(b\cR)$ such that $\cA^{1,\a}(b\cR)= (\ker A) \oplus \cB$
and $A$ maps $\cB$ isomorphically onto $\cC^{1,\a}_{\R}(b\cR)$. 
By the implicit function theorem in Banach spaces [Ca] the solutions 
of the equation $\Phi(\rho,f)=0$ in a small \nbd\ of 
$(\rho_0,f_0) \in \cC^3(P) \times \cA^{1,\a}(b\cR)$ 
are of the form $f(\rho,t)=f_0+ (t,\varphi(\rho,t))$, where
$t\in \ker A$ and $\varphi$ is a $\cC^1$-operator with image in
$\cB$ such that $\varphi(\rho_0,0)=0$ (and hence $f(\rho_0,0)=f_0$). 
Setting $t=0$ we obtain functions 
$f_\rho=f(\rho,0) \in\cA^{1,\a}(b\cR)$ for $\rho\in\cC^3(P)$
near $\rho_0$ such that $f_\rho$ depends differentiably on $\rho$
and satisfies $\rho(f_\rho,g_0)=0$ on $b\cR$. This concludes 
the proof of proposition 3.1.
\endpr

\beginsection 4. Families of inner functions 
on bordered Riemann surfaces.

In this section we prove theorem 1.4. The essential ingredient
is the following result which is possibly of independent interest. 
We use the notation established in section 1.
If $J$ is a complex structure of class $\cC^{k-1,\a}$ on $\cR$, 
we denote by $\cA^{k,\a}(\cR,J)$ the space of all $J$-holomorphic 
functions of order $\cC^{k,\a}$ (that is, their derivatives of 
order $k$ are H\"older continuous of order $\a$).

\proclaim 4.1 Proposition: Let $\cR$ be a smooth bordered 
surface of genus $g_{\cR}$ and with $m\ge 1$ boundary 
components. Fix $\a\in (0,1)$. Let $J_0\in \cJ_\cR^\a$ 
be a complex structure on $\cR$ and  let
$f_0\in \cA^{1,\a}(\cR,J_0)$ be an inner function 
on $(\cR,J_0)$ with multiplicity $\ge 2g_{\cR}+m-1$. 
Then for each $J\in \cJ_\cR^\a$ sufficiently close to 
$J_0$ there is an inner function $f_J\in  \cA^{1,\a}(\cR,J)$
near $f_0$, with $f_J$ depending continuously on $J$
and $f_{J_0}=f_0$.

\demo Remark: As already mentioned, Ahlfors [Ahl]
constructed inner functions of multiplicity $2g_{\cR}+m$
on any bordered Riemann surface. Proposition 4.1 shows 
that such functions are stable under small perturbations of 
the complex structure. On the other hand this need not be true
for the Ahlfors function $f_p$ which maximizes the derivative
at a given point $p\in \cR$ since the degree of $f_p$ may
depend on $p$. 
\endpr

Assuming proposition 4.1 for a moment we can prove
theorem 1.4 as follows. Fix a number $0<\a<1$ 
and let $F_0=(f_0,g_0)\colon \cR\to \bar U \times \C$ 
be an embedding of class $\cC^{1,\a}$ which is
$J_0$-holomorphic on $\IR$ for some complex structure
$J_0 \in \cJ_\cR^\a$. Thus $f_0$ is an inner function 
on $(\cR,J_0)$ of multiplicity $\ge 2g_\cR+m-1$. 
For each $J\in \cJ_\cR^\a$ sufficiently near $J_0$ 
proposition 4.1 provides an inner function $f_J$
on $(\cR,J)$ which is $\cC^1$-close to $f_0$.
We can also approximate $g_0$ in the $\cC^1$-sense
by $J$-holomorphic functions $g_J$ (this is trivial
since there is no boundary condition on $g_J$). 
If the approximations are sufficiently
close (which is the case when $J$ is close enough to $J_0$),
the $J$-holomorphic map 
$F_J=(f_J,g_J)\colon \cR\to \bar U\times \C$ 
is $\cC^1$-close to $F_0$ and hence is an embedding.
Hence the Riemann surface $(\cR,J)$ is of class $\cF$
for all $J$ sufficiently close to $J_0$. 
This proves theorem 1.4 provided that 
proposition 4.1 is correct.
\endpr

\demo Proof of proposition 4.1:
For each complex structure $J\in \cJ_\cR^\a$ we denote
by $\dibar_J$ the corresponding $\dibar$-operator which
maps $\cC^{1,\a}$-functions on $\cR$ to $(0,1)$-forms
of class $\cC^\a$ according to the formula
$$ 2\, \dibar_J(f) = df+i df\circ J. $$ 
Denote the space of such forms by $\Omega_{0,1}^\a(\cR,J)$. 
Consider the Banach manifolds 
$$ \eqalign{ \cW &= \{(f,J)\colon f \in \cC^{1,\a}(\cR),\ 
                     |f|=1{\rm\ on\ } b\cR,\ 
                     J\in \cJ_\cR^\a \},\cr
   \cW^{0,1} &=	\{(\omega,J)\colon\ \omega \in 
   \Omega_{0,1}^\a(\cR,J),\ J\in \cJ_\cR^\a \}. \cr}
$$
Let $\Phi\colon \cW\to \cW^{0,1}$ be the operator 
$\Phi(f,J)=(2\dibar_J f,J)$. The set 
$$ \cW^h =\{(f,J)\in \cW\colon \dibar_J(f)=0\}= 
          \{(f,J)\colon \Phi(f,J)=(0,J)\} 
$$
consists of $J$-holomorphic inner functions on $\cR$
for all complex structures $J$ on $\cR$. Denote by 
$\pi\colon \cW\to \cJ_\cR^\a$ the projection
onto the second factor.

We claim that $\Phi$ is a $\cC^1$-map of Banach manifolds
which is a submersion with finite corank at each point 
$(f,J)\in \cW^h$ for which $\k(f)\ge 2g_{\cR}+m-1$.
Once this is proved, the implicit function theorem [Ca] 
shows that $\cW^h$ is a Banach submanifold of $\cW$ 
in a \nbd\ of each such point $(f,J)$ and the projection 
$\pi \colon \cW^h \to \cJ_\cR^\a$ is locally near 
$(f,J)$ a trivial Banach fibration with finite dimensional 
fibers. The proposition then follows immediately since it amounts 
to choosing a local section of this fibration passing 
through $(f,J)$.

To find the derivative $D\Phi(f,J)(g,K)$ in the direction of
a tangent vector $(g,K)$ to $\cW$ at $(f,J)$ we choose
a local $\cC^1$ path $(f_t,J_t)\in \cW$ for $|t|<\e$,
with $f_0=f$, $J_0=J$, ${d\over dt}|_{t=0} f_t=g$ and 
${d\over dt}|_{t=0} J_t=K$. Differentiating the equations 
$|f_t|^2=1$ resp.\ $J_t^2=-Id$ with respect to $t$ at $t=0$ 
we see that ${\rm Re}(g \bar f)=0$ on $b\cR$ and $JK+KJ=0$.
The derivative of $\Phi$ equals
$$	\eqalign{ D\Phi(f,J)(g,K) &= {d\over dt}|_{t=0} \Phi(f_t,J_t) \cr
    &= \bigl( dg+idg J+idf  K,K\bigr) \cr
	&= \bigl( 2\dibar_Jg+idf K, K). \cr}
$$
From this formula we see immediately that $D\Phi(f,J)$ is continuous
in $(f,J)$ and hence $\Phi$ is of class $\cC^1$.
Moreover we see that $D\Phi(f,J)$ is surjective 
\iff\ any $\omega\in \Omega_{0,1}^\a(\cR,J)$
equals $\omega = 2\dibar_Jg + idf K$ for some 
$g\in \cC^{1,\a}(\cR)$ with ${\rm Re}(g \bar f)=0$ on $b\cR$.
Since $(idf K)J =-idf JK=dfK= -i(idf K)$ (here we used 
$KJ=-JK$ and $idfJ=-df$), the form $idfK$ is of type
$(0,1)$ with respect to $J$. Hence it suffices 
to see that $g\to\dibar_J g$ is surjective
as a map 
$$ \{g\in \cC^{1,\a}(\cR)\colon 
   {\rm Re}(g\bar f)=0\ {\rm on\ } b\cR \} 
   \to 	\dibar_J g \in \Omega_{0,1}^\a(\cR,J). \eqno(4.1)
$$
Surjectivity of this map at points $(f,J)\in \cW$ with
$\k(f)\ge 2g_{\cR}+m-1$ is guaranteed by the theorem in 
[Kop, p.\ 33] together with Corollary II in [Kop, p.\ 30]. 
(Essentially the result in [Kop] is that we can solve
any non-homogeneous Cauchy-Riemann equation on a fixed 
bordered Riemann surface $(\cR,J)$ subject to a Riemann-Hilbert 
boundary condition, provided that the associated index is 
sufficiently large, which in our case means that $f$ must
have multiplicity $\ge 2g+m-1$.)  This shows that
$\Phi$ is a submersion at such points as claimed. 
Furthermore, if $(f,J)\in \cW^h$, the kernel of 
the map (4.1) equals 
$$
	\{g\in  \cC^{1,\a}(\cR)\colon   \dibar_J g=0,\ 
   {\rm Re}(g\bar f)=0 {\rm \ on \ } b\cR\}.
$$
By the cited result in [Kop] this space 
has dimension $2\kappa(f)-(2g+m-2)$.
This completes the proof.
\endpr

\bigskip
\ni\bf References. \rm
\medskip

\ii{[Ahl]}  L.\ Ahlfors:
Open Riemann surfaces and extremal problems on compact subregions.
Comment.\ Math.\ Helv.\ {\bf 24} (1950), 100--134. 

\ii{[Ale]}  H.\ Alexander: 
Explicit imbedding of the (punctured) disc into $\C^2$. 
Comment.\ Math.\ Helv.\ {\bf 52} (1977), 539--544.

\ii{[AS]} 	M.\ F.\ Atiyah, I.\ M.\ Singer:
The index of elliptic operators. I. 
Ann.\ of Math.\ (2) {\bf 87} (1968), 484--530. 

\ii{[BS]} H.\ Behnke, F.\ Sommer:  Theorie der analytischen 
Funktionen einer komplexen Ver\"anderlichen, 3rd.\ ed.
Springer-Verlag, New York, 1965.

\ii{[Ca]} H.\ Cartan: Calcul diff\'{e}rentiel.
Hermann, Paris, 1967,

\ii{[EG]} Y.\ Eliashberg, M.\ Gromov: Embeddings of Stein manifolds.
Ann.\ Math.\ {\bf 136}, 123--135 (1992).

\ii{[\v CG]} M.\ \v Cerne, J.\ Globevnik:
On holomorphic embedding of planar domains into $\C^2$.
J.\ d'Analyse Math.\ {\bf 8} (2000), 269--282.

\ii{[FK]} H.\ M.\ Farkas, I.\ Kra: Riemann surfaces. 
Second ed. Graduate Texts in Mathematics, 71.
Springer, New York, 1992.

\ii{[FP1]} F.\ Forstneri\v c and J.\ Prezelj:
Oka's principle for holomorphic fiber bundles with sprays.
Math.\ Ann.\ {\bf 317} (2000), 117-154.

\ii{[FP2]} F.\ Forstneri\v c and J.\ Prezelj:
Oka's principle for holomorphic submersions with sprays.
Math.\ Ann.\ {\bf 322}(2002), 633-666.   

\ii{[FP3]} F.\ Forstneri\v c and J.\ Prezelj:
Extending holomorphic sections from complex subvarieties.
Math.\ Z.\ {\bf 236} (2001), 43--68.

\ii{[Glo]} J.\ Globevnik: On Fatou-Bieberbach domains.
Math.\ Z.\ {\bf 229} (1998), 91--106.

\ii{[GS]} J.\ Globevnik, B.\ Stens\o nes:
Holomorphic embeddings of planar domains into $\C^2$.
Math.\ Ann.\ {\bf 303} (1995), 579--597. 

\ii{[Gou]} T.\ Gouma: 
Ahlfors functions on non-planar Riemann surfaces whose double 
are hyperelliptic. 
J.\ Math.\ Soc.\ Japan {\bf 50} (1998), 685--695. 

\ii{[Gra]} H.\ Grauert:   
Holomorphe Funktionen mit Werten in komplexen Lieschen Gruppen. 
Math.\ Ann.\ {\bf 133} (1957), 450--472. 

\ii{[Gro1]} M.\ Gromov:
Pseudoholomorphic curves in symplectic manifolds. 
Invent.\ Math.\ {\bf 82} (1985), 307--347. 

\ii{[Gro2]} M.\ Gromov:
Oka's principle for holomorphic sections of elliptic bundles.
J.\ Amer.\ Math.\ Soc.\ {\bf 2}, 851-897 (1989).

\ii{[GR]} C.\ Gunning, H.\ Rossi:
Analytic functions of several complex variables.
Prentice-Hall, Englewood Cliffs, 1965.

\ii{[HT]}  D.\ C.\ Hill, D.\ Taiani: Families of analytic discs 
in $\C^n$ with boundaries in a prescribed CR manifold.
Ann.\ Sc.\ Norm.\ Sup.\ Pisa\ cl.\ sc.\ {\bf 5} (1978), 327--380.

\ii{[HLS]} H.\ Hofer, V.\ Lizan, J.-C.\ Sikorav:
On genericity for holomorphic curves 
in four-dimensional almost-complex manifolds.
J.\ Geom.\ Anal.\ {\bf 7} (1998), 149--159.

\ii{[H\"or]} L.\ H\"ormander:
An Introduction to Complex Analysis in Several Variables, 3rd ed.
North Holland, Amsterdam, 1990.

\ii{[Kop]} W.\ Koppelman: 
The Riemann-Hilbert problem for finite Riemann surfaces.
Comm.\ Pure Appl.\ Math.\ {\bf  12} (1959), 13--35.

\ii{[Lan]} S.\ Lang:\ Differential and Riemannian Manifolds.
Grad.\ Texts in Math.\ {\bf 160}, 
Springer, New York, 1995.

\ii{[Lau]} H.\ B.\ Laufer:
Imbedding annuli in $\bf C^2$. 
J.\ Analyse Math.\ {\bf  26} (1973), 187--215.

\ii{[NN]} A.\ Newlander, L.\ Nirenberg: 
Complex analytic coordinates in almost complex manifolds. 
Ann.\ of Math.\ (2) {\bf 65} (1957), 391--404. 

\ii{[Ru1]}  W.\ Rudin: 
Pairs of inner functions on finite Riemann surfaces. 
Trans.\ Amer.\ Math.\ Soc.\ {\bf 140} (1969), 423--434.

\ii{[Ru2]} W.\ Rudin: Functional Analysis.
McGraw-Hill, New York, 1973. 

\ii{[Sik]} J.\ C.\ Sikorav:
Proof that every torus with one hole can be properly
holomorphically embedded in $\C^2$.
Preprint, Oct.\ 1997 (unpublished).

\ii{[SS]} M.\ Schiffer, D.\ C.\ Spencer: 
Functionals on Finite Riemann Surfaces.
Princeton Univ.\ Press, Princeton, 1954.

\ii{[Sch]} J.\ Sch\"urmann:
Embeddings of Stein spaces into affine spaces 
of minimal dimension.
Math.\ Ann.\ {\bf 307}, 381--399 (1997).

\ii{[Spr]} G.\ Springer: Introduction to Riemann surfaces.
Addison-Wesley, Reading, Mass., 1957.

\ii{[Ste]} J.-L.\ Stehl\'e: Plongements du disque dans $\C^2$.
(S\'eminaire P.\ Lelong (Analyse), pp.\ 119--130), 
Lect.\ Notes in Math.\ {\bf 275}, Springer, Berlin--New York, 1970. 

\ii{[Stn]} B.\ Stens\o nes:
Fatou-Bieberbach domains with $\cC^\infty$-smooth boundary. 
Ann.\ of Math.\ (2) {\bf 145} (1997), 365--377.

\ii{[Tsa]} V.\ V.\ Tsanov:
On hyperelliptic Riemann surfaces and doubly generated function algebras. 
C.\ R.\ Acad.\ Bulgare Sci.\ {\bf 31} (1978), no.\ 10, 1249--1252.

%
%
%  Addresses
%
%
\bigskip\medskip
\itemitem{\it Address:} Institute of Mathematics, Physics and Mechanics,
University of Ljub\-ljana, Jadranska 19, 1000 Ljubljana, Slovenia

\bye